\documentclass{amsart}

\newtheorem{theorem}{Theorem}[section]

\newtheorem{corollary}[theorem]{Corollary}
\newtheorem{lemma}[theorem]{Lemma}

\theoremstyle{definition}
\newtheorem{definition}[theorem]{Definition}
\newtheorem{example}[theorem]{Example}

\theoremstyle{remark}
\newtheorem{remark}[theorem]{Remark}

\numberwithin{equation}{section}

\newcommand{\F}{\ensuremath{\mathcal{F}}}

\newcommand{\dank}{\textsf{Acknowledgments.\ }}

\begin{document}

%%% COMANDO ESPACO ENTRA AS LINHAS
%\baselineskip30pt

\title{On  polar foliations and fundamental group}

\author{Marcos M. Alexandrino}

%\thanks{ }

\address{ 
Instituto de Matem\'{a}tica e Estat\'{\i}stica\\
Universidade de S\~{a}o Paulo,\\ % \hfill\break\indent
 Rua do Mat\~{a}o 1010,05508 090 S\~{a}o Paulo, Brazil}
\email{marcosmalex@yahoo.de, malex@ime.usp.br}

\thanks{The  author was  supported by CNPq and partially supported by FAPESP. }


\date{ January 2011.}


\begin{abstract}
In this work we investigate the relation between the fundamental group of a complete Riemannian manifold $M$ and the quotient between the Weyl group and reflection group of a polar action on $M$, as well as the relation between the fundamental group of $M$  and the quotient between the lifted Weyl group and lifted reflection group. 
As  applications we give  alternative proofs of two results. The first one, due to  the  author and T\"{o}ben, implies that 
a polar action does not admit exceptional orbits, if $M$ is simply connected. The second result, due to Lytchak, implies that  the orbits  of a polar foliation are closed and embedded if $M$ is simply connected.
 All results are proved in the more general case of polar foliations. 
\end{abstract}

%\end{frontmatter}

\maketitle

\section{Introduction}
In this section, we  state our main results in Theorem \ref{theorem-surjectivehomeomorphism-polar}, Theorem 
\ref{theorem-surjectivehomeomorphism-desingularization},  Corollary \ref{corollary-AlexandrinoToben} and Corollary \ref{corollary-Lytchak}.

%A singular foliation $\F$ on a complete Riemannian manifold $M$ is 
%a \emph{singular Riemannian foliation} if each geodesic  perpendicular to one leaf   is perpendicular to every leaf it meets.

A singular Riemannian foliation $\F$ on a complete Riemannian manifold $M$ is called  a \emph{polar foliation}    
if for each regular point $p$, there is an immersed submanifold $\Sigma_{p}$, called \emph{section},  that passes through $p$ and that meets all the leaves and always perpendicularly. It follows that $\Sigma_p$ is totally geodesic and that the dimension of $\Sigma_{p}$ is equal to the codimension of the regular leaf $L_{p}$; for  definitions and properties of singular Riemannian foliations and polar foliations see \cite{Molino,AlexToeben2} and \cite{Alex2,Alex4,Toeben}.

By choosing an appropriate  section $\Sigma$, one can define a subgroup of isometries of $\Sigma$ called \emph{Weyl group} $W(\Sigma)$ of the section $\Sigma$ so that if $w\in W(\Sigma)$ then $w(x)\in L_{x}$. Therefore $W(\Sigma)$ describes how the leaves of $\F$ intersect the section $\Sigma$; see Section \ref{section-weyl-group}.  As expected, when $\F$ is  the partition of a Riemannian manifold $M$ into the orbits of a \emph{polar action}  $G\times M\to M$ (i.e., an isometric action with sections), the Weyl group $W(\Sigma)$ is the usual Weyl group $N/Z$, where $N=\{g\in G | \, g(x)\in \Sigma, \forall x\in \Sigma \}$ and $Z=\{g\in G | \, g(x)=x, \forall x\in \Sigma \}.$  The intersection of the singular leaves of $\F$ with the section $\Sigma$ is a union of totally geodesic hypersurfaces, called \emph{walls} and the reflections in the walls are elements of the Weyl group $W(\Sigma)$. The group generated by these reflections is called \emph{reflection group} $\Gamma(\Sigma).$ As explained in  Remark \ref{remark-diferencagruposreflexao}, the definition of $\Gamma(\Sigma)$ can be different from the usual definition of reflection group on a manifold. We are now able to state our first result.

\begin{theorem}
\label{theorem-surjectivehomeomorphism-polar}
Let $\F$ be a polar foliation on a complete Riemannian manifold $M$, $\Sigma$ a  section, $W=W(\Sigma)$ the Weyl group of $\Sigma$ and $\Gamma=\Gamma(\Sigma)$ the reflection  group of $\Sigma$. Then
there exists a surjective homomorphism $\rho: \pi_{1}(M)\to W/\Gamma.$ 
\end{theorem}

A direct consequence of Theorem \ref{theorem-surjectivehomeomorphism-polar} is that $W(\Sigma)=\Gamma(\Sigma)$ when $M$ is simply connected.
The fact that $W(\Sigma)=\Gamma(\Sigma)$ does not imply previous results about polar foliations on simply connected spaces; see details in Example \ref{subsection-example-so3}. In order to solve this problem, we lift the Weyl group $W(\Sigma)$ and reflection group $\Gamma(\Sigma)$ to subgroups of the isometries of the Riemannian universal cover of $\Sigma$, called \emph{lifted Weyl group} $\widetilde{W}$ and \emph{lifted reflection group} $\widetilde{\Gamma}$; see Definition \ref{def-lifted-Weyl-group} and Definition \ref{def-lifted-reflection-group}. 
This lead us to the second main result of this paper.

\begin{theorem}
\label{theorem-surjectivehomeomorphism-desingularization}
Let $\F$ be a polar foliation on a complete Riemannian manifold $M$. Let  $\widetilde{W}$ and $\widetilde{\Gamma}$ denote the lifted Weyl group and the lifted  reflection  group respectively.  Then 
there exists a surjective homomorphism $\widetilde{\rho}: \pi_{1}(M)\to \widetilde{W}/\widetilde{\Gamma}.$ 
\end{theorem}

Theorem \ref{theorem-surjectivehomeomorphism-desingularization} allow us to give   alternative proofs to the next two results, where the first one is due to the author and T\"{o}ben \cite{AlexToeben} and the second one is due to Lytchak \cite{Lytchak}.

\begin{corollary}
\label{corollary-AlexandrinoToben}
Let $\F$ be a polar foliation on a  simply connected complete Riemannian manifold. Then each regular leaf has trivial holonomy. 
\end{corollary}

\begin{corollary}
\label{corollary-Lytchak}
Let $\F$ be a polar foliation on a  simply connected complete Riemannian manifold. Then the leaves of $\F$ are closed and embedded.
\end{corollary}

This paper is organized as follows. In Section \ref{section-weyl-group} we recall the definition of Weyl group and present the definition of lifted Weyl group. 
In Sections \ref{section-proof-theoremA} we prove  Theorem \ref{theorem-surjectivehomeomorphism-polar}. In Section  \ref{section-proof-theoremB} we   prove  Theorem \ref{theorem-surjectivehomeomorphism-desingularization} using the construction of Section \ref{section-proof-theoremA}. 
In Section \ref{section-proof-of-the-corollaries}  we provide  alternative proofs of Corollaries \ref{corollary-AlexandrinoToben} and \ref{corollary-Lytchak}. In particular to read this section it suffices to read Section \ref{section-weyl-group} and accept   Theorem \ref{theorem-surjectivehomeomorphism-desingularization}. Finally, in Section \ref{section-Infinitesimally-polar}   we note that the above results admit  natural generalizations to the class of infinitesimally polar foliations.

\dank 
The author thanks  Dr.~Dirk T\"{o}ben and Professor Gudlaugur Thorbergsson for very helpful discussions. In particular the author is grateful to Prof.~Thorbergsson for his consistent support  and for suggesting one of the questions that motivates this work.

%%%%%%%%%%%%%%%%%%%%%%%
%%%%%%%%%%%% Weyl and lifted groups
%%%%%%%%%%%%%%%%%%%%%%%%%%%%%%%%%%%%%%%%%%%%%%%%%%%%%%%%%%%%%%%%%

\section{Weyl and lifted Weyl groups}
\label{section-weyl-group}
Let $\F$ be a polar foliation on a complete Riemannian manifold $M$.
It follows from \cite{Alex2} that, for each point $x_{\alpha},$ we can find a neighborhood $U_{\alpha}$ (an \emph{isoparametric neighborhood}) such that $\F|_{U_{\alpha}}$ is diffeomorphic to a polar foliation on an Euclidean space, i.e., an \emph{isoparametric foliation}. Therefore we can find a totally geodesic submanifold $\sigma_{\alpha}\subset U_{\alpha}$ transverse to the plaques that we  call \emph{local section}.
 
Now  consider $\{U_{\alpha}\}$ an open covering of $M$ by isoparametric neighborhoods.
If $U_{\alpha}\cap U_{\beta}\neq \emptyset$, then we  find neighborhoods $V^{\alpha}_{\alpha}\subset \sigma_{\alpha}$, 
$V^{\alpha}_{\beta}\subset \sigma_{\beta}$ and an isometry $\varphi_{\beta,\alpha}:V^{\alpha}_{\alpha}\cap C_{\alpha}\to V^{\alpha}_{\beta} \cap C_{\beta}$
where $C_{\alpha}$ (respectively $C_{\beta}$) is a \emph{Weyl chamber} of $\sigma_{\alpha}$ (respectively $\sigma_{\beta}$), i.e., 
$C_i\subset \sigma_i$  is
the closure in $\sigma_{i}$ of a connected component of the set of nonsingular points in $\sigma_{i},$ for $i=\alpha,\beta$. 
Since $C_{\alpha},$ $C_{\beta}$ are Weyl chambers, we  extend the map $\varphi_{\alpha,\beta}$ to an isometry 
$\varphi_{\beta,\alpha} :V^{\alpha}_{\alpha}\to V^{\alpha}_{\beta}$. 
If $\beta$ is a curve contained in a regular leaf, we can cover $\beta$ by open neighborhoods $U_{i}$
and set $\varphi_{[\beta]}:=\varphi_{n,n-1}\circ\ldots\circ\varphi_{2,1}$. The germ of $\varphi_{[\beta]}$
is then called \emph{holohomy map} and depends only on the homotopy class $[\beta]$ of $\beta$.
The source (domain) of $\varphi_{[\beta]}$ can contain singular points.

We  define the  Weyl pseudogroup  of a section $\Sigma$ as the pseudogroup generated by the holonomy maps whose source and target are contained in the section $\Sigma$. By the appropriate choice of section, this pseudogroup turns out to be a group that is called  \emph{Weyl group} $W(\Sigma)$; see \cite{Toeben}.
From now on we fix this section where the Weyl group is well defined as the section $\Sigma$.
By definition, if $w\in W(\Sigma)$ then $w(x)\in L_{x}$ and $W(\Sigma)$ describes how the leaves of $\F$ intersect the section $\Sigma$. As expected, when $\F$ is  the partition of a Riemannian manifold $M$ into the orbits of a polar action $G\times M\to M$, the Weyl group $W(\Sigma)$ is the usual Weyl group $N/Z$, where $N=\{g\in G | \, g(x)\in \Sigma, \forall x\in \Sigma \}$ and $Z=\{g\in G | \, g(x)=x, \forall x\in \Sigma \}$; see \cite{PTlivro}. 

It follows from \cite{Alex2} that the intersection of the singular leaves of $\F$ with the section $\Sigma$ is a union of totally geodesic hypersurfaces, called \emph{walls} and that the reflections in the walls are elements of the Weyl group $W(\Sigma)$.
We  define the \emph{reflection group} $\Gamma(\Sigma)$  as the group generated by these reflections. Again when $\F$ is the partition of a Riemannian manifold into the orbits of a polar action, $\Gamma(\Sigma)$ is the usual reflection group in the section.

\begin{remark}
\label{remark-diferencagruposreflexao}
Note that, in the definition of a reflection of a polar foliation, we do not assume that fixed point set $\Sigma_{r}$ of a reflection $r$ separates $\Sigma$, as usual in the  literature. Also note that what we call \emph{wall} is a  hypersurface contained in the intersection of the singular leaves of $\F$ with $\Sigma$. Therefore a wall is a connected component of the fixed point set of a reflection $r$, but there can exist a connected component of the fixed point set of $r$ that is not a wall. This is also different from the usual definition of walls of reflections. These definitions  will coincide when we consider the lifted reflection group.
\end{remark}

 Consider $w\in W(\Sigma)$. In what follows we will define a isometry $\tilde{w}:\widetilde{\Sigma}\to\widetilde{\Sigma}$ on the Riemannian universal cover of $\Sigma$ such that $\pi_{\Sigma}\circ\tilde{w}=w\circ\pi_{\Sigma}$ where $\pi_{\Sigma}:\widetilde{\Sigma}\to \Sigma$ is the Riemannian covering map.

\begin{definition}
Let $p$ be  regular point of the section $\Sigma$  and $\tilde{p}$ a point of $\widetilde{\Sigma}$ such that $\pi_{\Sigma}(\tilde{p})=p$.
Consider  a curve $c$ in $\Sigma$ that joins $p$ to $w(p)$ and $\delta$ be a curve in $\Sigma$ with $\delta(0)=p$. Let  $c *w\circ\delta$ be the concatenation  of $c$ and the curve $w\circ\delta$, $\widetilde{\delta}$ the lift of  $\delta$ starting at $\tilde{p}$ and $\widetilde{( c * w \circ \delta)}$ the lift of $c *w\circ\delta$ starting at $\tilde{p}$ .
Then we  set
\begin{equation}
\label{equation-new-def-tildew}
\tilde{w}_{[c]}(\widetilde{\delta}(1)):=\widetilde{( c * w \circ \delta)}(1).
\end{equation}
The isometry $\tilde{w}_{[c]}$ defined in equation \eqref{equation-new-def-tildew} is  called \emph{a lift of $w$ along $c$}. Clearly it depends on the homotopy class $[c]$.
\end{definition}

\begin{definition}
\label{def-lifted-Weyl-group}
We  call \emph{lifted Weyl group} $\widetilde{W}$ the group of isometies on the universal covering $\widetilde{\Sigma}$ generated by all isometries $\tilde{w}_{[c]}$ constructed above. In other words,
$\widetilde{W}=<\tilde{\omega}_{[c]}> $ for all $ w\in W$ and for all curves $ c:[0,1]\to \Sigma,$ such that $c(0)=p$ and $c(1)=w(p).$
\end{definition}

\begin{remark}
Sometimes we will need to lift an isometry $w$ that is  naturally characterized by its action on a neighborhood of a point $x\in\Sigma$ far from the fixed point $p$. In this case, it is convenient to consider the following procedure. Let  $c_{0}$ be a curve in $\Sigma$ that joins $p$ to $x$ and $c_{1}$ a curve in $\Sigma$ that joins $x$ to $w(x)$.
Let $\delta$ be a curve in $\Sigma$ with $\delta(0)=x$, $c_{0} * \delta$ the concatenation  of $c_0$ and $\delta$ and 
$\widetilde{(c_{0}*\delta)}$ the lift of  $c_{0}*\delta$ starting at $\tilde{p}$.
Then we  define 
\begin{equation}
\label{equation-def-tildew}
\hat{w}(\widetilde{(c_{0}*\delta)}(1)):=\widetilde{(c_{0}* c_{1}* w \circ \delta)}(1).
\end{equation}
Clearly the construction of $\tilde{w}$ depends on the homotopy class $[c_0]$ and $[c_{1}].$
It is not difficult to check that the group generated by all these isometries is the lifted Weyl group $\widetilde{W}$. 
\end{remark}

\begin{remark}
\label{remark-about-lifted-Weyl-group}
The lifted Weyl group was used by T\"{o}ben \cite{Toeben} in his study of desingularization of polar foliations.  One can prove that each element $\tilde{w}\in\widetilde{W}$ is a lift of some isometry $w\in W$.

\end{remark}

We are now able to construct the appropriate reflection group on $\widetilde{\Sigma}$.

\begin{definition}
\label{def-lifted-reflection-group}
Let $H$ be a wall in $\Sigma$. Consider a reflection $r\in\Gamma(\Sigma)$ in this wall $H$ and a point $x$ near to $H$. Let $c_1$ be  a curve joining $x$ to $w(x)$ in a simply connected neighborhood. Now define $\tilde{r}$ by equation (\ref{equation-def-tildew}). The isometry $\tilde{r}$ is then a reflection in a hypersurface $\widetilde{H}$ such that $\pi_{\Sigma}(\widetilde{H})=H$. The group generated by all reflections $\tilde{r}$ is called
\emph{lifted reflection group} $\widetilde{\Gamma}$. 
\end{definition}

\begin{remark}
\label{remark-lifted-reflection-group}
One can prove that  $\widetilde{\Gamma}$ is a normal subgroup of $\widetilde{W}$. 
From the definition we have that $\pi_{1}(\Sigma)$ is a subgroup of $\widetilde{W}$. 
If $W(\Sigma)$ is the Weyl group of $\Sigma$ then one can check that  $\pi_{1}(\Sigma)$ is a normal subgroup and $W(\Sigma)=\widetilde{W}/\pi_{1}(\Sigma)$. Nevertheless, in general $\pi_{1}(\Sigma)$ is not entirely contained in $\widetilde{\Gamma}$. A simple example is again  the action of $SO(3)$ on itself by conjugation (see Example \ref{subsection-example-so3}). In particular $\widetilde{W}/\widetilde{\Gamma}$ does not need to be isomorphic to $W/\Gamma$. 

\end{remark}

\begin{remark}
\label{remark-widegamma-classico}
It follows from Davis \cite[Lemma 1.1]{Davis} that $\widetilde{\Gamma}$ is a \emph{reflection group in the classical sense}, i.e., is a discrete group generated by reflections acting properly, effectively on $\widetilde{\Sigma}$. In particular a Weyl chamber $C$ of $\widetilde{\Gamma}$ is a fundamental domain of the action of $\widetilde{\Gamma}$; see \cite[Theorem 4.1]{Davis}.
\end{remark}

%%%%%%%%%%%
%%%%%%PROOF THEOREM A
%%%%%%%%%%%%%%%%%%%%%%%%%%%%%%%%%%%%
\section{Proof of Theorem \ref{theorem-surjectivehomeomorphism-polar}}
\label{section-proof-theoremA}

\subsection{The construction of the map $\rho$}
\label{subsection-construction-rho}
In this subsection we construct a map $\rho$ that associate each loop $\alpha$ with an element of $W/\Gamma$. 
Let $p$ be a regular point and $\alpha:[0,1]\to M$ be a loop with $\alpha(0)=p=\alpha(1)$.

We start by considering the following objects:

\begin{enumerate}
\item A partition $0=t_{0}<t_{1}<\ldots<t_{n}=1$ and an open covering $\{U_{i}\}$ of $\alpha$ by \emph{isoparametric neighborhood}  (i.e., $\F|_{U_{i}}$ is diffeomorphic to an isoparametric foliation) such that $\alpha_{i}:=\alpha|_{I_{i}=[a_{i},b_{i}]}\subset U_{i}$, where $a_{i}=t_{i-1}$ and $b_{i}=t_{i}$;
\item  $\sigma_i\subset U_i$  \emph{a local section}, i.e., a normal neighborhood of a section $\Sigma_i$, and $C_i\subset \sigma_i$ a \emph{Weyl chamber}, i.e.,
the closure in $\sigma_{i}$ of a connected component of the set of nonsingular points in $\sigma_{i};$ 
\item $\pi_{i}:U_{i}\to C_{i}$ the \emph{chamber map}, i.e., the continuous map such that $L_{y}\cap C_{i}=\{\pi_{i}(y)\}$. 

\end{enumerate}

Now we consider an holonomy map
$\varphi_{i}:=\varphi_{[\beta_{i}]}:V_{i}^{i}\subset \sigma_{i}\to V_{i+1}^{i}\subset \sigma_{i+1}$ with the following properties:
\begin{enumerate}
\item $\beta_{i}$ is a curve in an isoparametric neighborhood $U_{i,i+1}\supset U_{i}\cup U_{i+1},$
\item $\pi_{i}\circ\alpha_{i}(b_{i})\in V_{i}^{i}$ and $\pi_{i+1}\circ\alpha_{i+1}(a_{i+1})\in V_{i+1}^{i},$
\item $\varphi_{i}\circ\pi_{i}\circ\alpha_{i}(b_{i})=\pi_{i+1}\circ\alpha_{i+1}(a_{i+1}).$
\end{enumerate}

We  also consider holonomies maps $\psi_{i}:B_{i}\subset \Sigma\to \sigma_{i}\subset \Sigma_{i}$ with $\psi_{1}=Id=\psi_{n+1}$.

Finally we define, for $1\leq i\leq n,$ the following objects:
\begin{enumerate}
\item $\gamma_{i}:=\psi_{i}^{-1}\circ\pi_{i}\circ\alpha_{i}$,
\item  $ w_{i}:=\psi_{i+1}^{-1}\circ\varphi_{i}\circ\psi_{i}.$ 
\end{enumerate}

We have then constructed a \emph{$W$-loop} $(\gamma_{i},w_{i})$ \emph{with base point} $p\in\Sigma,$ i.e.,    
\begin{enumerate}
\item a sequence $0=t_{0}<\cdots<t_{n}=1,$    
\item continuous paths $\gamma_{i}:[t_{i-1},t_{i}]\rightarrow \Sigma,$ $1\leq i\leq n,$  
\item elements $w_{i}\in W$ defined in a neighborhood of $\gamma_{i}(t_{i})$ for $1\leq i \leq n$ such that $\gamma_{1}(0)=w_{n}\gamma_{n}(1)=p$ and $w_{i}\gamma_{i}(t_{i})=\gamma_{i+1}(t_{i}),$ where $1\leq i\leq n-1.$
\end{enumerate}

Note that we can use the $W$-loop $(\gamma_{i},w_{i})$ to construct a connected curve $\gamma$ with endpoint $w_{1}^{-1}\cdots w_{n}^{-1}p$. For example, if we have
a $W$-loop $(\gamma_1,w_1)$, $(\gamma_2,w_2)$, $(\gamma_{3},w_{3})$, we have the connected curve 
$\gamma_{1}* (w_{1}^{-1}\gamma_{2})*( w_{1}^{-1}w_{2}^{-1}\gamma_{3})$ with end point $w_{1}^{-1}w_{2}^{-1}\gamma_{3}(1)=w_{1}^{-1}w_{2}^{-1}w_{3}^{-1}p$. 
We finally define the map $\rho$ as follows:
\begin{equation}
\label{eq-definition-rho}
\rho(\alpha):=(w_{1})^{-1}\cdot (w_{2})^{-1}\cdots (w_{n})^{-1}\Gamma.
\end{equation}

%%%%%%%%%%%%%%%%%%%%%%
%%%%%%%%%%%%%%%%%%%%%%%%%%%%SECAO 3-SUBSECAO weel defined
%%%%%%%%%%%%%%%%%%%%%%%%%%%% 

\subsection{The map $\rho:\pi_{1}(M,p)\to W/\Gamma$ is well defined and surjective}
\label{subsection-rho-welldefined}
We start by checking that the definition of the map  $\rho$  in equation (\ref{eq-definition-rho}) does not depend on the choice of the Weyl chamber $\{C_i\}$, the local sections $\{\sigma_i\}$ and the maps $\{\psi_i\}$.

For another choice of local sections $\{\hat{\sigma}_{i}\},$ Weyl chambers $\{\hat{C}_{i}\}$ and maps $\{\hat{\psi}_{i}\}$ we have, as explained in the last subsection,
a $W$-loop $(\hat{\gamma}_{i},\hat{w}_{i})$.

By straightforward calculations on can check the next equation
\begin{equation}
\label{eq-hatwi-wi}
\hat{w}_{i}=(g_{i+1})^{-1}(\mathrm{r}_{i+1}^{i})^{-1}w_{i}\mathrm{r}_{i}^{i}g_{i}, 
\end{equation}
where $g_1=e=g_{n+1}$, $g_i\in W$, $\mathrm{r}_{j}^{i}\in\Gamma$.
Equation \eqref{eq-hatwi-wi} then implies
$$(w_{1})^{-1}\cdot (w_{2})^{-1}\cdots (w_{n})^{-1}\Gamma=(\hat{w}_{1})^{-1}\cdot (\hat{w}_{2})^{-1}\cdots (\hat{w}_{n})^{-1}\Gamma.$$

Note that the construction does not depend on the open covering $\{U_{i}\}$. In fact, for another open covering $\{\widehat{U}_{j}\}$ we  assume, without lost of generality, that for each $j$ there exists $i$ such that $\widehat{U}_{j}\subset U_{i}$.
For this new open covering we have a $W$-loop $(\hat{\gamma}_{j},\hat{w}_{j}).$
Now one can infer equation (\ref{eq-hatwi-wi}) by considering the appropriate  
\emph{subdivision} of  the $W$-loop $(\gamma_{i},w_{i})$, i.e., by adding new points to the interval $[0,1],$ by taking the restriction of the $\gamma_{i}$ to these new intervals and $w=id$ at the new points.
Finally we have to verify that the construction does not depend on the homotopy class of $\alpha$. Consider a fixed choice of $\{U_{i}\}$, $\{\sigma_{i}\}$, $\{C_{i}\}$ and $\{\psi_{i}\}$. Then one can check that for each curve $\alpha^{s}$ near to $\alpha$ we produce a $W$-loop $(\gamma_{i}^{s},w_{i}^{s})$ where $w_{i}^{s}=w_{i}$. 
Therefore $\rho(\alpha^{s})=\rho(\alpha).$

We have proved that the map  $\rho:\pi_{1}(M,p)\to W/\Gamma$ is well defined. By choosing the appropriate open covering, one can see that $\rho$ is also a homomorphism of groups. We conclude this section proving that it is surjective.

Consider $w\in W$ and let $\gamma$ be a curve that joins $p$ to $w(p)$.  
Now consider $\beta\subset L_{p}$ a curve such that joins $p$ to $w(p)$ and such that $\varphi_{[\beta]}=w.$ Then consider the concatenation $\alpha=\gamma*\beta^{-1}$. With the appropriate choice of $\{U_{i}\}$, $\{\sigma_{i}\}$, $\{C_{i}\}$ and $\{\psi_{i}\}$, we produce the $W$-loop $(\gamma,w^{-1})$ and hence  $\rho(\alpha)=(w^{-1})^{-1}\Gamma=w\Gamma.$

%%%%%%%%%%%%%%%%%%
%%%%%%%%%%%%%%%%%%
%%%%%% SECAO 4 PROVA DO TEOREMA B
%%%%%%%%%%%%%%%%%%

\section{Proof of Theorem \ref{theorem-surjectivehomeomorphism-desingularization}}

\label{section-proof-theoremB}

In this section we will use the notations and definitions of the previous sections.

\subsection{The construction of the map $\tilde{\rho}$ and an example}
\label{subsection-definicaofacil-tilderho}

Let $p$ be a regular point and $\alpha:[0,1]\to M$ be a loop with $\alpha(0)=p=\alpha(1)$. 
In this subsection we associate  $\alpha$ with  an element $\tilde{\rho}(\alpha)\in \widetilde{W}/\widetilde{\Gamma},$  
 briefly sketch the rest of the proof of Theorem \ref{theorem-surjectivehomeomorphism-desingularization} and illustrate the lifted groups and maps $\rho$ and $\tilde{\rho}$ in an example.

Consider $\{U_{i}\}$,$\{\sigma_{i}\}$, $\{C_{i}\}$ and $\{\psi_{i}\}$ defined in  Section \ref{section-proof-theoremA}.
We will also assume that, for $1\leq i<n,$ the open set $B_{i}\cup B_{i+1}$ is contained in a simply connected neighborhood $B_{i,i+1}$ where $B_{i}=\psi^{-1}_{i}(\sigma_{i})$. As explained in Section \ref{section-proof-theoremA}, we  produce a
$W$-loop $(\gamma_{i},w_{i})$ and with this $W$-loop we  construct a connected curve 
$\gamma$ such that $w\gamma(0)=\gamma(1)$, where $w=(w_{1})^{-1}\cdot (w_{2})^{-1}\cdots (w_{n})^{-1}$. Finally we  set
\begin{equation}
\label{eq-def-tilderho-facil}
\tilde{\rho}(\alpha):=\tilde{w}_{[\gamma]}\tilde{\Gamma},
\end{equation}
where $\tilde{w}_{[\gamma]}$ is the lift of the isometry $w$ along the curve $\gamma$; recall  equation \eqref{equation-new-def-tildew}. 
Different choices of  Weyl chambers will produce different curves $\gamma$ and isometries $w$, but the fact that we are considering the quotient space $\widetilde{W}/\widetilde{\Gamma}$ will imply that $\tilde{\rho}$ is well defined; see Subsection \ref{subsection-tilderho-welldefined}.

To prove that $\tilde{\rho}$ is surjective, for each $\tilde{w}\in \widetilde{W}$ we have to find a loop $\alpha$ such that  $\tilde{\rho}(\alpha)=\tilde{w}\widetilde{\Gamma}$. Here we have to use the fact that each $\tilde{w}$ is a lift $\tilde{w}_{[c]}$ of an isometry $w$. Then we can choose $\alpha=c *\beta^{-1}$ where  $\varphi_{[\beta]}=w$.

\begin{example}
\label{subsection-example-so3}
Assuming that $\tilde{\rho}$ is well defined, we now illustrate the lifted groups and maps $\rho$ and $\tilde{\rho}$ in 
 the simple example of the action of $G=SO(3)$ on itself by conjugation. In this example the section $\Sigma$ is isomorphic to $S^{1}$. The identity $Id\in SO(3)$ is the only singular point in $\Sigma$ and the antipodal point $p\in\Sigma$ is the only point such that the orbit $G(p)$ is exceptional. Note that $W=\Gamma=\{e, w\}$. Therefore \emph{the fact that $W(\Sigma)=\Gamma(\Sigma)$ does not imply that all regular orbits have trivial holonomy.}
Clearly $\widetilde{\Sigma}=\mathbb{R}$ and we can identify the elements of $\widetilde{\Gamma}$ with the reflection in integer numbers.
Note that the group $\pi_{1}(\Sigma)$ of deck transformations is not entirely contained in $\widetilde{\Gamma}$, because the translation of length 1 is a deck transformation and is not contained in $\widetilde{\Gamma}$.  We can identify $\tilde{p},$ the lift of the exceptional point $p$, with the point $1/2\in\mathbb{R}$. 
Note that, by the definition of lifted isometries, $R_{1/2}$, the reflection in $1/2$, is an element of $\widetilde{W}$. In fact all elements of $\widetilde{W}$ can be generated by $R_{1/2}$ and $R_{1}$, the reflection in the point 1. 

Now we exam the maps $\rho$ and $\tilde{\rho}$. Consider $\beta$ a curve in $G(p)$ such that $\varphi_{[\beta]}=w$. Set $\alpha_{1}=\beta^{-1}$.  Then by the process explained in Section \ref{section-proof-theoremA} we produce a $W$-loop $(p,w^{-1})$. Therefore $\rho(\alpha_{1})=w\Gamma=\Gamma.$ 
One can check that $\tilde{\rho}(\alpha_{1})=R_{1/2}\widetilde{\Gamma}$.
Finally set $\alpha_{2}=\gamma*\beta^{-1}$ where $\gamma$ is a parametrization of $\Sigma=S^{1}$. Then by the process explained in Section \ref{section-proof-theoremA} we produce a 
$W$-loop $(\gamma, w^{-1})$. Therefore $\rho(\alpha_{2})=w\Gamma=\Gamma.$ One can check that $\tilde{\rho}(\alpha_{2})=R_{1}\widetilde{\Gamma}=\widetilde{\Gamma}$.
\end{example}

%%%%%%%%%%%%%%%%%%%%%
%%%%%%%%%%% SECAO 
%%%%%%%%%%%%%%%%%%%%

\subsection{The map $\tilde{\rho}:\pi_{1}(M,p)\to \widetilde{W}/\widetilde{\Gamma}$ is well defined}
\label{subsection-tilderho-welldefined}

Let  $(\gamma_{i},w_{i})$ be the $W$-loop defined in Subsection \ref{subsection-definicaofacil-tilderho}.
In order to prove that $\tilde{\rho}$ is well defined, we will need an equivalent definition of $\tilde{\rho}$ written in terms of  lifts of the isometries $w_{i}$.

Let $\tilde{\gamma}_{1}$ be the lift of $\gamma_{1}$ starting at  $\tilde{p}$. We  join $\gamma_{1}(b_{1})$ to $\gamma_{2}(a_{2})$ by a curve $c_{1}\subset B_{1,2}$ and define $\tilde{\gamma}_{2}$ the lift of $\gamma_{2}$ starting at $\widetilde{(\gamma_{1}* c_{1})}(1)$.
We define the lift $\tilde{w}_{1}$ replacing the curve $c_{0}$ with the curve $\gamma_{1}$ in equation \eqref{equation-def-tildew}. Since $B_{1,2}$ is simply connected, the definition of $\tilde{\gamma}_{2}$ and $\tilde{w}_{1}$ does not depend on the choice of $c_{1}$. We  join $\gamma_{i}(b_{i})$ with $\gamma_{i+1}(a_{i+1})$ by a curve $c_{i}\subset B_{i,i+1}$ and, following the same procedure,  define by induction
$\tilde{\gamma}_{1},\ldots,\tilde{\gamma}_{n}$ and $\tilde{w}_{1},\ldots,\tilde{w}_{n-1}$. The fact that $B_{i,i+1}$ is simply connected imply that this construction does not depend on the choice of $c_{i}$. In order to define $\tilde{w}_{n}$ set $c_{0}:=\gamma_{1}* c_{1}\cdots * c_{n-1}*\gamma_{n}$, set $c_{n}:=c_{0}^{-1}$ and replace $c_{1}$ with $c_{n}$ in equation \eqref{equation-def-tildew}.
Note that, for another choice of curves $c_1,\ldots, c_{n-1}$, we would have an homotopic curve joining $p$ to $\gamma_{n}(b_{n})$. Therefore the definition of $\tilde{w}_{n}$ does not depend on the choice of the curves $c_1,\ldots, c_{n-1}$. 
The next lemma is not difficult to check.

\begin{lemma}
\label{lemma-definition-tilderho}
Let $\tilde{w}_i$ defined above. Then
$\tilde{\rho}(\alpha)=(\tilde{w}_{1})^{-1}\cdot (\tilde{w}_{2})^{-1}\cdots (\tilde{w}_{n})^{-1}\widetilde{\Gamma}.$
\end{lemma}

Consider another choice of local sections $\{\hat{\sigma}_{i}\}$, Weyl chambers $\{\hat{C}_{i}\}$ and holonomies  $\{\hat{\psi}_{i}\}$. As before, we
 also assume that, for $1\leq i<n,$ the open set $\hat{B}_{i}\cup \hat{B}_{i+1}$ is contained in a simply connected neighborhood $\hat{B}_{i,i+1}$ where $\hat{B}_{i}=\hat{\psi}_{i}^{-1}(\hat{\sigma}_{i})$. As explained in Section \ref{section-proof-theoremA}, we can produce another 
$W$-loop $(\hat{\gamma}_{i},\hat{w}_{i})$ and then we have a $\widetilde{W}$-loop $(\widetilde{\hat{\gamma}}_{i},\widetilde{\hat{w}}_{i}).$
Using equation \eqref{equation-def-tildew} one can prove that 
\begin{equation}
\label{eq-tildehatwi-tildewi-1}
 \tilde{g}_{i+1}\widetilde{\hat{w}}_{i}= (\tilde{\mathrm{r}}_{i+i}^{i})^{-1}\tilde{w}_{i}\tilde{\mathrm{r}}_{i}^{i} \tilde{g}_{i}. 
\end{equation}
where $\tilde{\mathrm{r}}_{j}^{i}\in\widetilde{\Gamma}$, $\tilde{g}_{i}\in \widetilde{W}$ and $\tilde{g}_1=Id=\tilde{g}_{n+1}$.
This equation and the fact that $\widetilde{\Gamma}$ is a normal subgroup of $\widetilde{W}$ then imply  
\begin{equation}
\label{eq-tilde-rho-independence}
(\tilde{w}_{1})^{-1}\cdot (\tilde{w}_{2})^{-1}\cdots (\tilde{w}_{n})^{-1}\tilde{\Gamma}=(\widetilde{\hat{w}}_{1})^{-1}\cdot (\widetilde{\hat{w}}_{2})^{-1}\cdots (\widetilde{\hat{w}}_{n})^{-1}\tilde{\Gamma}.
\end{equation}

From Lemma \ref{lemma-definition-tilderho} we conclude that the construction does not depend on the choice of $\{\sigma_{i}\}$, $\{C_{i}\}$ and $\{\psi_{i}\}$.  

By subdivision of $W$-loop we see that the definition also does not depend on the covering $\{U_{i}\}$. We have concluded that the map $\tilde{\rho}$ is well defined. The independence of  the homotopy class of $\alpha$ is again proved considering a fixed choice of $\{U_{i}\}$, $\{\sigma_{i}\}$, $\{C_{i}\}$ and $\{\psi_{i}\}$ and noting that, for each curve $\alpha^{s}$ near to $\alpha$, we produce an $W$-loop $(\gamma_{i}^{s},w_{i}^{s})$ where $w_{i}^{s}=w_{i}$. 
Therefore $\tilde{\rho}(\alpha^{s})=\tilde{\rho}(\alpha).$

%\begin{remark}
%Note that the assumption that the Weyl pseudogroup of $\Sigma$ is a Weyl group $W(\Sigma)$ is not necessarily  in the definition of the Weyl and reflection %lifted groups. Also note that the proof of Theorem \ref{theorem-surjectivehomeomorphism-desingularization} works, with small modifications, without this %assumption.    
%\end{remark}

%%%%%%%%%%%%%%%%%%%%%%%%%%%%%%%%%%%%%%%%%%%%%%%%%%5
%%%%%%%%%%%%%%%%%%%%%%%%%% SECAO 5 PROVA DOS COROLARIOS
%%%%%%%%%%%%%%%%%%%%%%%%%%%%%%%%%%%%%%%%%%%%%%%%%%%%%%%% 

\section{Proof of the corollaries}
\label{section-proof-of-the-corollaries}

%%%%%%%%%%%%%%%%%%%%%%%%%%%%%%%%%%%%%%%%%%%
\subsection{Proof of Corollary \ref{corollary-AlexandrinoToben}}

Assume that there exists a regular leaf with a nontrivial holonomy, i.e., assume that there exists a regular point $p\in\Sigma$ and $w\in W$ such that $w(p)=p$ and $w\neq e$.

On one hand we  define a lift of $w$ as  $\tilde{w}(\widetilde{\delta}(1)):=\widetilde{w\circ\delta}(1).$ Note that $\tilde{w}(\tilde{p})=\tilde{p}$ and $\tilde{w}$ is not the identity.

On the other hand, since $M$ is simply connected, Theorem \ref{theorem-surjectivehomeomorphism-desingularization} implies that $\widetilde{W}=\widetilde{\Gamma}$.  
It follows from Remark \ref{remark-widegamma-classico} that the only points that can be fixed by elements of $\widetilde{\Gamma}$ are singular points, i.e., points whose projection by $\pi_{\Sigma}$ are singular points of $\Sigma$.
Therefore $\tilde{w}(\tilde{p})\neq \tilde{p}$ and we have arrived at a contradiction. 

%%%%%%%%%%%%%%%%%%
\subsection{Proof of Corollary \ref{corollary-Lytchak}}

Due to equifocality  \cite{AlexToeben2}, in order to prove Corollary \ref{corollary-Lytchak}, it  suffices to prove that the leaves are embedded.

Assume that $L_{p}$ is not embedded. Then we can find a sequence $\{w_{n}\}\subset W$ such that $w_{n}(p)\to p$. Define 
$\tilde{w}_{n}(\widetilde{\delta}(1)):=\widetilde{(c_{n}* w_{n}\circ\delta)}(1)$, where $c_{n}$ are curves with small lengths contained in a simply connected neighborhood of $p$. Then one can check that $\tilde{w}_{n}(\tilde{p})\to \tilde{p}$. 

On the other hand, since $M$ is simply connected,  Theorem \ref{theorem-surjectivehomeomorphism-desingularization} implies that $\tilde{W}=\tilde{\Gamma}$.  
By Remark \ref{remark-widegamma-classico},   the sequence $\{\tilde{w}_{n}(\tilde{p})\}$ can not converge to $\tilde{p}$.  Hence we have arrived at a contradiction.

%%%%%%%%%%%%%%%%%%%%%%%%%%%%%%%
%%%%%%% SECAO 6 INFINITESIMALLY POLAR FOLIATIONS
%%%%%%%%%%%%%%%%%%%%%%%%%%%%%%%%%%%%%%%

\section{Infinitesimally polar foliations}
\label{section-Infinitesimally-polar}

A singular Riemannian foliation $\F$ is called \emph{infinitesimally polar foliation} if the restriction  of $\F$  to each slice is diffeomorphic (by the composition of the exponential map with a linear map) to  an isoparametric foliation.
Lytchak and Thorbergsson \cite{LytchakThorbergsson1} proved that $\F$ is an infinitesimally polar foliation if and only if for each point $x\in M$ we can find a neighborhood $U$ of $x$ such the leave space of the restricted foliation $\F|_{U}$ is an orbifold. Since the leave space of an isoparametric foliation is a Coxeter orbifold, they concluded that the leave space $U/\F|_{U}$ is an orbifold if and only if it is a Coxeter orbifold.
Typical examples of infinitesimally polar foliations are polar foliations (see \cite{Alex2}), singular Riemannian  foliations without horizontal conjugate points (see \cite{LytchakThorbergsson1,LytchakThorbergsson2}) and singular Riemannian foliations with codimension lower then three (see \cite{LytchakThorbergsson2}). 

\begin{remark}
Let $\F$ be an infinitesimally polar foliation on a complete Riemannian manifold $M$ and $\{U_{\alpha}\}$ a covering of $M$ by isoparametric 
 neighborhoods. We can find  submanifolds 
$\sigma_{\alpha}$ transverse to the plaques that we still call \emph{local sections}.
Using the same argument as in  \cite[Subsection 3.2]{LytchakThorbergsson2} we can change the
metric of local sections $\sigma_{\alpha}$ so that the holonomy maps turn out to be isometries.
Like in the classical theory of foliations, the structure of $M/\F$ is described by the \emph{pseudogroup of holonomies} of $\F$, i.e., 
the pseudogroup generated by holonomy maps acting on the disjoint union of local sections.

Now assume that $M/\F$ is equivalent to $\Sigma/W$  where $\Sigma$ is 
a connected Riemannian manifold and $W$ is a subgroup of isometries of $\Sigma$;  
 for details about pseudogroup see Salem \cite[Appendix D]{Molino}.
Then $W$ is called \emph{Weyl group} of $\F$. Let $\{\psi_{\alpha}\}$ be a maximal collection
of isometries of open subsets of $\Sigma$ onto local sections of $\F$ that gives the equivalence
between $\Sigma/W$ and $M/\F$.  Let $r$  be a reflection in a wall of a local section $\sigma_{\alpha}$. 
Then $\mathrm{r}:=\psi_{\alpha}^{-1}\circ r\circ\psi_{\alpha}\in W$  is called a \emph{reflection}
 on $\Sigma$ and $\Gamma$ is the group generated by these reflections.
One can also define the \emph{lifted Weyl group} $\widetilde{W}$ and the 
\emph{lifted reflection group} $\widetilde{\Gamma}$ as 
in Definitions \ref{def-lifted-Weyl-group} and \ref{def-lifted-reflection-group}. 
The same proofs of  Theorem \ref{theorem-surjectivehomeomorphism-polar}
and Theorem \ref{theorem-surjectivehomeomorphism-desingularization} allow us to conclude that
\emph{there exist surjective homomorphisms $\rho: \pi_{1}(M)\to W/\Gamma$ and  $\widetilde{\rho}: \pi_{1}(M)\to \widetilde{W}/\widetilde{\Gamma}.$}
Also the same proof of Corollary \ref{corollary-Lytchak} implies that \emph{the leaves of $\F$ are closed and embedded and $M/\F$ is a good orbifold, if $M$ is simply connected}. 
\end{remark}

Taking the above remark into account,  one can adapt  the proof of Corollary \ref{corollary-AlexandrinoToben} and  get  the next result due to Lytchak \cite{Lytchak}. 
\begin{corollary}
\label{habilitation-lytchak}
Let  $\F$ be a singular Riemannian foliation on a simply connected complete Riemannian manifold $M$. Assume that $M/\F$ is a good orbifold $\Sigma/W$.
Then  the regular leaves have trivial holonomy.
\end{corollary}

\begin{remark}
It is not difficult to see that, if the regular leaves have trivial holonomies, then $M/\F$ is a Coxeter orbifold; see \cite{AlexToeben}. Therefore, if $M$ is simply connected, $M/\F$ is a Coxeter good orbifold. 
This result is  used by the author and Javaloyes \cite{AlexMiguel} to prove that  $M/\F$ admits nontrivial closed geodesic if $M/\F$ is a compact orbifold and $\F$ is a singular Riemannian foliation on a simply connected complete manifold $M$.
\end{remark}

%%%%%%%%%%%%%%%%%%%%%%%%%%%%%%%%%%%%%%%%%%%%%%%%

\bibliographystyle{amsplain}

\end{document}